\documentclass[a4paper,11pt]{article}
\usepackage{graphicx}
\usepackage{multirow}
\usepackage{color}
\usepackage{latexsym}
\usepackage{amssymb}
\usepackage{amsfonts}
\usepackage{amsmath}
\usepackage{amsthm}
\usepackage{indentfirst}
\usepackage{mathrsfs}
\usepackage{tikz}
\usepackage{comment}
\usetikzlibrary{patterns}
\usetikzlibrary{hobby,decorations.pathreplacing,arrows.meta}
\usepackage{tikz-qtree,tikz-qtree-compat}

\makeatletter
\def\underbrace#1{%
   \@ifnextchar_{\tikz@@underbrace{#1}}{\tikz@@underbrace{#1}_{}}}
\def\tikz@@underbrace#1_#2{%
   \tikz[baseline=(a.base)] {\node[inner sep=2] (a) {\(#1\)};
   \draw[line cap=round,decorate,decoration={brace,amplitude=4pt}]
     (a.south east) -- node[pos=0.5,below,inner sep=7pt] {\(\scriptstyle #2\)} (a.south west);}}
 
\def\overbrace#1{%
   \@ifnextchar^{\tikz@@overbrace{#1}}{\tikz@@overbrace{#1}^{}}}
\def\tikz@@overbrace#1^#2{%
   \tikz[baseline=(a.base)] {\node[inner sep=2] (a) {\(#1\)};
   \draw[line cap=round,decorate,decoration={brace,amplitude=4pt}]
     (a.north west) -- node[pos=0.5,above,inner sep=7pt] {\(\scriptstyle #2\)} (a.north east);}}
\makeatother

\usepackage{fullpage}

\newtheorem {theorem}{Theorem}[section]

\newtheorem {proposition}{Proposition}[section]
\newtheorem {corollary}{Corollary}[section]

\newtheorem {example}{Example}[section]

\newcommand{\N}{\mathbb{N}}

\title {A generating tree with a single label for permutations avoiding the vincular pattern $1-32-4$}
\author{ Matteo Cervetti
\thanks{Dipartimento di Matematica, Università di Trento.
{\tt \ matteo.cervetti@unitn.it}}}
 
\date{} 
 
\begin{document}

\maketitle

\begin{abstract}
In this paper we continue the study  of permutations avoiding the vincular pattern $1-32-4$ by constructing a generating tree with a single label for these permutations. This construction finally provides a clearer explanation of why a certain recursive formula found by Callan actually counts these permutations, insofar as this formula was originally obtained only as a consequence of a very intricated bijection with a certain class of ordered rooted trees. This responds to a theoretical issue already raised by Duchi, Guerrini and Rinaldi.  As a byproduct, we also obtain an algorithm to generate all these permutations and we refine their enumeration according to a simple statistic, which is the number of right-to-left maxima to the right of 1.  
\end{abstract}

\section{Introduction}

The study of patterns in permutations has grown in the last two decades to one of the most active trends of research in contemporary combinatorics.  The study of permutations which are constrained by not having one or more subsequences ordered in various prescribed ways has been historically motivated by the problem of sorting permutations by means of certain devices. However, nowadays research on this topic is further fueled by its intrinsic combinatorial difficulty and from the plentiful appearances of patterns in several very different disciplines, such as algebra, geometry, analysis, computer science, mathematical physics and computational biology, and many others. For this reason, it is reasonable to believe that this field of research will continue growing for a long time to come.

  More recently, the notion of vincular (or generalized) pattern in permutations has been  considered. Whereas an occurrence of a classical pattern $\pi$ in a permutation $\sigma$ is simply a (non-necessarily consecutive) subsequence of $\sigma$ whose items are in the same relative order as those in $\pi$, in an occurrence of a vincular pattern, some items of that subsequence may be required to be adjacent in the permutation. For instance, the classical pattern $1234$ simply corresponds to an increasing subsequence of length 4, whereas an occurrence of the generalized pattern $1-23-4$ would require the middle two letters of that sequence to be adjacent in $\sigma$. Thus, the permutation $23145$ contains $1234$ but not $1-23-4$. The study of vincular patterns provides significant additions to the extensive literature on classical patterns. In fact, the non-classical vincular patterns are likely to provide richer connections to other combinatorial structures than the classical ones do. Other than combinatorics, vincular patterns find applications in other scientific topics such as, for instance, in the genome rearrangement problem, which is one of the major trends in bioinformatics and biomathematics. 
  
The main issue we adress in this paper concerns permutations avoiding a well known vincular pattern of length $4$, namely the pattern $1-32-4$. The paper is organized as follows. 

In Section $\ref{VP}$ we fix the main notations and terminology which we will use throughout the paper. 

 In Section $\ref{1-32-4avoiders}$ we describe our contribution to the study of permutations avoiding $1-32-4$. Permutations avoiding this pattern were actually enumerated by Callan in $\cite{C}$, who provided a recursive formula to count them only as a collateral consequence of a very intricated bijection with a certain class of ordered rooted trees. The problem of finding a more explicit justification for the recursive formula obtained by Callan has received some attention in the past decade (see $\cite{DGR}$).  In Section $\ref{1-32-4avoiders}$ we construct a generating tree with a single label for permutations avoiding the vincular pattern $1-32-4$, finally providing such a justification. 
 
 In order to introduct this topic, in Sections $\ref{PW}$ and $\ref{ECO}$ we collect the latest results about vincular pattern avoidance and we quickly review the ECO method, which is the framework needed to describe the construction in Section $\ref{1-32-4avoiders}$.

\section{The notion of vincular pattern}\label{VP}

 Let $n\in \mathbb{N}^{*}$ and denote $[n]=\{1,2,...,n\}$. For our purposes, a \emph{permutation of length} $n$ will be just a word with  $n$ distinct letters from the alphabet $[n]$. The length of a permutation $\sigma$ will be denoted by $|\sigma|$.  We will denote by $\mathcal{S}_{n}$ the set of all permutations of length $n$ and we will set $\mathcal{S}=\bigcup_{n\in \mathbb{N}}\mathcal{S}_{n}$ (where by convention we set $\mathcal{S}_{0}=\{\varepsilon\}$ and $|\varepsilon|=0$).
 
We begin by quickly recalling the classical definition of pattern in a permutation. In general, given any poset $\mathcal{P}$ and two non-empty words $\alpha$ and $\beta$ with length $k$ in the alphabet $\mathcal{P}$, we say that $\alpha$ and $\beta$ are \emph{order isomorphic}, and we will write $\alpha\sim \beta$ when, for every $i,j\in [k]$, $\alpha_{i}\leq \alpha_{j}$ if and only if $\beta_{i}\leq \beta_{j}$.  Let now $\sigma,\tau\in \mathcal{S}$ and suppose $i\in [|\tau|]^{|\sigma|}$. We say that $i$ is an $\emph{occurrence}$ of $\sigma$ in $\tau$ when $i_{1}<i_{2}<...<i_{|\sigma|}$ and $\tau_{i_{1}}\tau_{i_{2}}...\tau_{i_{|\sigma|}}\sim\sigma$. We say that $\sigma$ \emph{occurs as pattern} in $\tau$, and we write $\sigma\leq \tau$, when either $\sigma=\varepsilon$ or one can find an occurrence of $\sigma$ in $\tau$, we say that $\tau$ $\emph{avoids}$ $\sigma$ \emph{as a} $\emph{pattern}$ otherwise. It is routine to check that $\leq$ is a partial order relation, turning $\mathcal{S}$ into a poset, which is called the \emph{permutation pattern poset}. 

It is worth noting that the items in an occurrence of a classical pattern in a permutation are not required to be necessarily consecutive, forcing this further condition we obtain the consecutive pattern poset on $\mathcal{S}$. With the notations as above, we say that $i$ is a $\emph{consecutive}$ $\emph{occurrence}$ of $\sigma$ in $\tau$ when $i$ is a occurrence of $\sigma$ in $\tau$ (in the classical sense) and either $|\sigma|=1$ or $i_{j+1}=i_{j}+1$ for every $j\in [|\sigma|-1]$. In the same fashion, we say that $\sigma$ \emph{occurs as consecutive pattern} in $\tau$ when one can find a consecutive occurrence of $\sigma$ in $\tau$. The resulting poset on $\mathcal{S}$ will be called the $\emph{consecutive}$ $\emph{pattern}$ $\emph{poset}$. Usually, the consecutive pattern poset reveals a much simpler structure than the classical pattern poset. For instance,  the M\"obius function of the consecutive permutation pattern poset is completely understood, whereas it is largely unknown in the classical case. 

Consecutive and classical patterns are special (actually, extremal) cases of the more general
notion of vincular patterns. Vincular patterns were introduced by Babson and Steingrimsson (under the name of generalized patterns), and constitute a vast intermediate continent between the two lands of consecutive patterns and classical patterns. Investigation of this intermediate notion could hopefully shed some light on differences and analogies between the extremal cases of consecutive and classical patterns.  An occurrence of a vincular pattern is basically an occurrence of that pattern in which entries are subject to given adjacency conditions. For a more formal definition of vincular pattern we follow $\cite{BF}$. 

A $\emph{dashed}$ $\emph{permutation}$ is a permutation in which some dashes are possibly inserted between any two consecutive letters. For instance, $5-13-42$ is a dashed permutation (of length 5). The type of a dashed permutation $\sigma$ such that $|\sigma|\geq 2$  is the $\{0, 1\}-$vector $r = (r_{1},...,r_{|\sigma|-1})\in \{0,1\}^{|\sigma|-1}$ such that $r_{i} = 0$ whenever there is no dash between $\sigma_{i}$ and $\sigma_{i+1}$ and $r_{i} = 1$ whenever there is a dash between $\sigma_{i}$ and $\sigma_{i+1}$, for every $i\in [|\sigma|-1]$. For example, the above dashed permutation $5-13-42$ has type $(1,0,1,0)$.  Given a dashed permutation $\sigma$, the $\emph{underlying}$ $\emph{permutation}$ of $\sigma$ is the permutation obtained by removing the dashes from $\sigma$. For instance, the compresed permutation of $5-13-42$ is the permutation $51342$. Let $\sigma$ be a dashed permutation. With the same notations as before, given some $i=(i_{1},...,i_{|\sigma|})\in [|\tau|]^{|\sigma|}$, we say that $i$ is an $\emph{occurrence}$ of $\sigma$ in $\tau$ when it is a (classical) occurrence of the underlying  permutation of $\sigma$ and $i_{j+1}=i_{j}+1$ whenever $j\in [|\sigma|-1]$ and $\sigma_{j}$ and $\sigma_{j+1}$ are not separated by a dash. We say that $\sigma$ \emph{occurs as a vincular pattern} in $\tau$ when one can find an occurrence of $\sigma$ in $\tau$, we say that $\tau$ \emph{avoids} $\sigma$ \emph{as a vincular pattern} otherwise.  Given a dashed permutation $\sigma$, the set of all permutations avoiding $\sigma$ will be denoted by $\mathcal{S}(\sigma)$ and the set of all $\sigma\in \mathcal{S}(\sigma)$ such that $|\sigma|=n$ will be denoted by $\mathcal{S}_{n}(\sigma)$.

\section{Previous work on vincular pattern avoidance}\label{PW}
  
In this section we provide a quick overview of the latest results in enumeration of permutations avoiding a vincular pattern. We refer to $\cite{Ki1}$ and $\cite{Ste}$ for a more detailed survey on this topic. The first systematic study of vincular patterns of length 3 was done by Claesson in $\cite{Cl}$ and permutations avoiding $\pi$ have been enumerated for every vincular pattern $\pi$ of length  3.  

 It is worth mentioning that the fact that, as it turns out, $\mathcal{S}(1-23)$ is counted by the Bell numbers shows that the analogous of the Stanley-Wilf conjecture does not hold for some vincular patterns. Moreover, the same fact shows that the conjecture of Noonan and Zeilberger stated in $\cite{NZ}$ is also false for vincular patterns, namely, the number of permutations avoiding a vincular pattern is not necessarily polynomially recursive.  
  
  As for vincular patterns of length 4, it is stated in \cite{Ste} that there are 48 symmetry classes of vincular patterns of length 4, and computer experiments show that there are at least 24 Wilf-equivalent classes (although their exact number seems to be unknown). For vincular non-classical patterns enumerative results are known for seven Wilf classes (out of at least 24) which are as follows:

\begin{itemize}
\item Elizalde and Noy $\cite{EN}$ gave the exponential generating functions for the number of occurrences of a consecutive pattern of length $4$ for three out of the seven Wilf-equivalence classes for consecutive patterns, namely the classes with representatives 1234, 1243 and 1342.
\item  Kitaev $\cite{Ki2}$ and Elizalde $\cite{E}$ decomposed the class $\mathcal{S}(\sigma -k)$ in a suitable  boxed product, where $\sigma$ is any consecutive pattern and $k=|\sigma|+1$. This decomposition allows to provide an expression for the exponential generating function of $\mathcal{S}(\sigma -k)$ in terms of the exponential generating function of $\mathcal{S}(\sigma)$. In particular, if $\sigma$ is any consecutive pattern of length 3, this, together with the results of Elizalde and Noy $\cite{EN}$, yields an explicit formula for the exponential generaing function $\mathcal{S}(\sigma-4)$, where $\sigma$ is any consecutive pattern of length 3. Since there are precisely two Wilf-equivalence classes of consecutive patterns of length 3, with representatives $123$ and $132$, the result of Kitaev yields explicit formulas for the exponential generating functions of  $\mathcal{S}(\pi)$ where $\pi$ is any vincular pattern Wilf-equivalent to $123-4$ or $132-4$. Explicitly, these formulas are
$$\exp\left( \frac{\sqrt{3}}{2}\int_{0}^{x}\frac{e^{\frac{t}{2}}}{\cos\left(\frac{\sqrt{3}}{2}t+\frac{\pi}{6}\right)}\right) \ \ \mathrm{and}\ \ \exp\left(\int_{0}^{x}\frac{dt}{1-\int_{0}^{t}e^{-\frac{u^{2}}{2}}du}\right)$$
for a vincular pattern Wilf-equivalent to $123-4$ or $132-4$, respectively.
\item Callan gave a recursion for $a_{n}=|\mathcal{S}_{n}(31-4-2)|$, which goes as follows. Set $a_{0}=c_{1}=1$ and 
\begin{itemize}
\item[1.] $a_{n}=\sum_{i=0}^{n-1}a_{i}c_{n-i}$ for $n\geq 1$.
\item[2.] $c_{n}=\sum_{i=0}^{n-1}ia_{(n-1),i}$ for $n\geq 2$.
\item[3.] $a_{n,k}=\begin{cases}\sum_{i=0}^{k}a_{i}\sum_{j=k-i}^{n-1-i}a_{(n-1-i),j} & 1\leq k\leq n-1\\ a_{n-1} & k=n\end{cases}$
\end{itemize}
\item Finally, Callan also showed in $\cite{C}$ that $|\mathcal{S}_{n}(1-32-4)|=\sum_{k=1}^{n}u(n,k)$ where, for every $1\leq k\leq n$, the triangle $u(n,k)$ satisfies the recurrence relation \begin{equation}\label{Callan}
u(n,k)=u(n-1,k-1)+k\sum_{j=k}^{n-1}u(n-1,j)
\end{equation}
with initial conditions $u(0,0)=1$ and $u(n,0)=0$ for every $n\geq 1$.
\end{itemize}

As far as the author knows, these seem to be the only explicit enumerative results concerning vincular patterns of length 4. An alternative, and we believe more explicative, proof of the recursive formula in Equation $(\ref{Callan})$ will be the main issue of the next sections.

\section{ECO method}\label{ECO}

In this section we provide a quick digression on ECO method, which proved to be the most suitable framework for describing the construction in Section $\ref{1-32-4avoiders}$.  The ECO method (Enumerating Combinatorial Objects method) was introduced by Barcucci, Del Lungo, Pergola and Pinzani $\cite{BDLPP}$ and it is quite a natural approach to generation and enumeration of combinatorial classes of objects according to their size. The main idea of the ECO method consists of looking for a way to grow objects from smaller to larger ones by making some local expansions, where each object should be obtained from a unique father so that this construction gives rise to a tree that allows us to recursively generate all the objects in the class. If the shape of this tree can be described with a simple rule there is hope for exact enumeration results, translating this description into equations for the generating function of the class. 

In the following we borrow standard terminology and notation about combinatorial classes from $\cite{FS}$.  Let $\mathcal{A}$ be a combinatorial class and denote by $\mathcal{A}_{n}$ the set of all $x\in \mathcal{A}$ such that $|x|_{\mathcal{A}}=n$ for every $n\in \mathbb{N}$. Let $\vartheta:\mathcal{A}\longrightarrow 2^{\mathcal{A}}$ be a map. We say that $\vartheta$ is an $\emph{ECO}$ $\emph{operator}$ on $\mathcal{A}$ when for every $n\in \mathbb{N}$ the following conditions hold:
\begin{itemize}
\item[(i)] $\vartheta(\mathcal{A}_{n})\subseteq 2^{\mathcal{A}_{n+1}}$.
\item[(ii)] for every $y\in \mathcal{A}_{n+1}$ one can find some $x\in \mathcal{A}_{n}$ such that $y\in \vartheta(x)$.
\item[(iii)] $\vartheta(x_{1})\cap \vartheta(x_{2})=\emptyset$ for every $x_{1},x_{2}\in \mathcal{A}_{n}$ such that $x_{1}\neq x_{2}$. 
\end{itemize}

In particular $\{\vartheta(x):\ x\in \mathcal{A}_{n}\}$ is a partition of $\mathcal{A}_{n+1}$. In other words, the map $\vartheta$ generates all the objects of the class $\mathcal{A}$ in such a way that each object in $\mathcal{A}_{n+1}$ is obtained from a unique object in $\mathcal{A}_{n}$ for every $n\in \mathbb{N}$. Actually, we can also characterize ECO operators as follows. We say that a map $\rho:\mathcal{A}\smallsetminus\mathcal{A}_{0}\longrightarrow\mathcal{A}$ is a $\emph{reduction}$ $\emph{operator}$ on $\mathcal{A}$ when $\rho(\mathcal{A}_{n+1})\subseteq \mathcal{A}_{n}$ for every $n\in \mathbb{N}$. The next proposition states that defining an ECO operator is essentially equivalent to defining a reduction operator. In the following, for a function $f:A\longrightarrow B$, we will denote by $f^{\longleftarrow}$ the function $f^{\longleftarrow}:B\longrightarrow 2^{A}$ assigning to each $b\in B$ its preimage $f^{\longleftarrow}(b)$ under $f$, i.e. the set of all $a\in A$ such that $f(a)=b$. 

\begin{proposition} Let $\mathcal{A}$ be a combinatorial class and $\vartheta:\mathcal{A}\longrightarrow 2^{\mathcal{A}}$ be a map. Then $\vartheta$ is an ECO operator on $\mathcal{A}$ if and only if $\vartheta=\rho^{\longleftarrow}$ for some reduction operator $\rho$ on $\mathcal{A}$. 
\end{proposition}

\proof Suppose $\vartheta$ is an ECO operator and $y\in \mathcal{A}\smallsetminus \mathcal{A}_{0}$, then $y\in \mathcal{A}_{n+1}$ for some $n\in \mathcal{N}$, hence $y\in \vartheta(x)$ for a unique $x\in \mathcal{A}_{n}$, because $\vartheta$ is an ECO operator, and we set $\rho(y)=x$. Let now $\rho:\mathcal{A}\longrightarrow\mathcal{A}$  denote the map assigning the object $\rho(y)$ to each $y\in \mathcal{A}$.  Then by definition $\rho(\mathcal{A}_{n+1})\subseteq \mathcal{A}_{n}$ for every $n\in \mathbb{N}$, furthermore, for every $x\in \mathcal{A}$, again by definition $y\in\vartheta(x)$ if and only if $x=\rho(y)$, equivalently $y\in \rho^{\longleftarrow}(x)$, hence $\vartheta(x)=\rho^{\longleftarrow}(x)$ and $\vartheta=\rho^{\longleftarrow}$. Conversely, suppose $\vartheta=\rho^{\longleftarrow}$ for some reduction operator $\rho$ on $\mathcal{A}$.  Take $n\in \mathbb{N}$, $x\in \mathcal{A}_{n}$ and $y\in \vartheta(x)=\rho^{\longleftarrow}(x)$, then $\rho(y)=x$, hence $y\in \mathcal{A}_{n+1}$ because $\rho(\mathcal{A}_{n+1})\subseteq \mathcal{A}_{n}$. Pick now $y\in \mathcal{A}_{n+1}$ and let $x=\rho(y)\in \mathcal{A}_{n}$, then by assumption $y\in \rho^{\longleftarrow}(x)=\vartheta(x)$. Finally, suppose $x_{1},x_{2}\in \mathcal{A}_{n}$ and $y\in \vartheta(x_{1})\cap\vartheta(x_{2})=\rho^{\longleftarrow}(x_{1})\cap \rho^{\longleftarrow}(x_{2})$, then $x_{1}=\rho(y)=x_{2}$. This proves that $\vartheta$ is an ECO operator on $\mathcal{A}$.  
\endproof

 Let $\vartheta$ be an ECO operator on $\mathcal{A}$. Say that $\mathcal{A}$ is $\emph{rooted}$ when $\mathcal{A}$ contains a unique object with minimum size. Suppose $\mathcal{A}$ is rooted. In this case, one can represent $\vartheta$ by means of a tree, called the $\emph{generating}$ $\emph{tree}$ of $\vartheta$ and denoted by $\mathcal{T}_{\vartheta}$, which is a rooted tree  having the objects of $\mathcal{A}$ as nodes, the object in $\mathcal{A}$ with minimum size as root and such that $\vartheta(x)$ is the set of children of $x$ for every $x\in \mathcal{A}$. This representation of $\vartheta$ can be useful for enumeration purpose when $\mathcal{T}_{\vartheta}$ displays enough regularity to be described by a so-called succession rule. Let $S$ be a set, let $a$ be an element of $S$, let $e$ be a sequence of maps in $S^{S}$ and let $p:S\longrightarrow\mathbb{N}$ be a map. The triple $(a,e,p)$ will be called the $\emph{succession}$ $\emph{rule}$ with $\emph{axiom}$ $a$, $\emph{production}$ $\emph{rule}$ $e$ and $\emph{production}$ $\emph{parameter}$ $p$ and it will be denoted by the symbol
$$\begin{cases}(a)\\ (k) \leadsto (e_{1}(k))(e_{2}(k))...(e_{p(k)}(k))\end{cases}$$
Let now $\mathcal{T}$ be a rooted tree with set of nodes $T$ and root $R$, let $\Omega=(a,e,p)$ be a succession rule and let $\ell:T\longrightarrow S$ be a map. We say that $\ell$ is a $\Omega-\emph{labelling}$ of $\mathcal{T}$ when $\ell(R)=a$ and, if $v\in T$, then $v$ has $p(\ell(v))$ children $v_{1},...,v_{p(\ell(v))}$ and $\ell(v_{i})=e_{i}(\ell(v))$ for every $i\in [p(\ell(v))]$. We say that $\vartheta$ and $\mathcal{A}$ are $\emph{described}$ by $\Omega$ when one can find some $\Omega-$labelling of $\mathcal{T}_{\vartheta}$. Suppose now $S\subseteq \mathbb{N}$ and $\ell$ is an $\Omega-$labelling of $\mathcal{T}_{\vartheta}$. Denote by $\mathcal{A}(z,u)$ the generating function
$$\mathcal{A}(z,u)=\sum_{\sigma\in \mathcal{A}}z^{|\sigma|}u^{\ell(\sigma)}=\sum_{n,k\geq 0}|\mathcal{A}_{n,k}|z^{n}u^{k}$$
where $\mathcal{A}_{n,k}$ denotes the set of all $\sigma\in \mathcal{A}_{n}$ such that $\ell(\sigma)=k$ for every $(n,k)\in \mathbb{N}^{2}$. We can translate the $\Omega-$description of $\mathcal{A}$ into a functional equation for $\mathcal{A}(z,u)$ as follows. Denote by $L_{\Omega}$ the unique $\mathbb{Z}[[z]]-$linear map $L_{\Omega}:\mathbb{Z}[[z,u]]\longrightarrow \mathbb{Z}[[z,u]]$ such that 
$$L_{\Omega}(u^{k})=\begin{cases}
u^{a} & k=0\\ u^{e_{1}(k)}+...+u^{e_{p(k)}(k)} & k\geq 1
\end{cases}$$ Then it is easily seen that 
$$\mathcal{A}(z,u)=\mathcal{A}(0,u)+zL_{\Omega}(\mathcal{A}(z,u)).$$
In the luckiest cases, these kind of equations can be solved using kernel type methods or other standard tools.

\section{Permutations avoiding the pattern $1-32-4$}\label{1-32-4avoiders}

The vincular pattern $1-32-4$ is a dashed version of the classical pattern $1324$, which attracted great attention among combinatorialists, as enumeration of permutations avoiding this classical pattern has proven to be one of the hardest open problems in permutation pattern combinatorics. Hopefully, our insight into the class of permutations avoiding the classical pattern could benefit from a closer study of permutations avoiding one of its vincular counterparts. However, these two  patterns seem also to display quite a different behaviour, for instance it is not difficult to see that $1-32-4$ is actually Wilf-equivalent to  $1-23-4$ (see $\cite{E}$), whereas this does not hold for their classical counterparts. 

Permutations avoiding $1-32-4$ are counted by sequence A113227 in $\cite{S}$, whose first ten terms are 1, 1, 2, 6, 23, 105, 549, 3207, 20577, 143239. As mentioned in Section $\ref{PW}$, an efficient bivariate recursive formula  to count permutations avoiding $1-32-4$ was first discovered by Callan in $\cite{C}$. This formula actually relies on a very intricated bijection (involving several contrived discrete structures defined ad hoc in order to break this transformation into somewhat simpler steps) between permutations of length $n$ avoiding $1-32-4$ and increasing ordered rooted trees on $n+1$ nodes with increasing leaves for every $n\geq 1$. An ordered rooted tree on $n+1$ nodes $\{0,1,2,...,n\}$ is $\emph{increasing}$ when every node is smaller than each of its children. If in addition its leaves are increasing from left to right we say that it has $\emph{increasing}$ $\emph{leaves}$. For instance, the figure below shows two increasing ordered rooted trees,
the first has increasing leaves while the second does not.
 
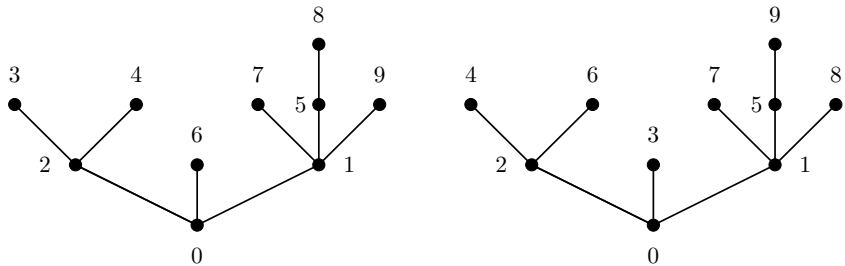
\begin{center}
\begin{tikzpicture} 
 \node at (0,0) {
\scalebox{0.8}{\begin{tikzpicture}
\draw [fill] (1,2) circle [radius=0.1];
\draw [fill] (2,1) circle [radius=0.1];
\draw [fill] (3,2) circle [radius=0.1];
\draw [fill] (4,0) circle [radius=0.1];
\draw [fill] (4,1) circle [radius=0.1];
\draw [fill] (5,2) circle [radius=0.1];
\draw [fill] (6,1) circle [radius=0.1];
\draw [fill] (6,2) circle [radius=0.1];
\draw [fill] (6,3) circle [radius=0.1];
\draw [fill] (7,2) circle [radius=0.1];
\node at (1,2.5) {3};
\node at (1.5,1) {2};
\node at (3,2.5) {4};
\node at (4,-0.5) {0};
\node at (4,1.5) {6};
\node at (6.5,1) {1};
\node at (5,2.5) {7};
\node at (6,3.5) {8};
\node at (7,2.5) {9};
\node at (5.7,2) {5};
\draw[thick] (4,0)--(2,1)--(1,2);
\draw[thick] (4,0)--(2,1)--(3,2);
\draw[thick] (4,0)--(4,1);
\draw[thick] (4,0)--(6,1)--(6,2)--(6,3);
\draw[thick] (6,1)--(5,2);
\draw[thick] (6,1)--(7,2);
\end{tikzpicture}}
 };
 
\node at (6,0) {
\scalebox{0.8}{\begin{tikzpicture}
\draw [fill] (1,2) circle [radius=0.1];
\draw [fill] (2,1) circle [radius=0.1];
\draw [fill] (3,2) circle [radius=0.1];
\draw [fill] (4,0) circle [radius=0.1];
\draw [fill] (4,1) circle [radius=0.1];
\draw [fill] (5,2) circle [radius=0.1];
\draw [fill] (6,1) circle [radius=0.1];
\draw [fill] (6,2) circle [radius=0.1];
\draw [fill] (6,3) circle [radius=0.1];
\draw [fill] (7,2) circle [radius=0.1];
\node at (1,2.5) {4};
\node at (1.5,1) {2};
\node at (3,2.5) {6};
\node at (4,-0.5) {0};
\node at (4,1.5) {3};
\node at (6.5,1) {1};
\node at (5,2.5) {7};
\node at (6,3.5) {9};
\node at (7,2.5) {8};
\node at (5.7,2) {5};
\draw[thick] (4,0)--(2,1)--(1,2);
\draw[thick] (4,0)--(2,1)--(3,2);
\draw[thick] (4,0)--(4,1);
\draw[thick] (4,0)--(6,1)--(6,2)--(6,3);
\draw[thick] (6,1)--(5,2);
\draw[thick] (6,1)--(7,2);
\end{tikzpicture}}
};
\end{tikzpicture} 
\end{center}

   Let $\mathcal{I}$ denote the combinatorial class of such trees. If, for every $1\leq k\leq n$,  we denote by $u(n)$ the number of trees in $\mathcal{I}_{n}$ and by $u(n,k)$ the number of trees in $\mathcal{I}_{n}$ such that the root has $k$ children, so that $u(n)=\sum_{k=1}^{n}u(n,k)$, then it is easily proved in $\cite{C}$ that the triangle $u(n,k)$ satisfies the recurrence relation 
 \begin{equation}\label{(1-32-4)-recursion}
 u(n,k)=u(n-1,k-1)+k\sum_{j=k}^{n-1}u(n-1,j)
 \end{equation}
when $1\leq k\leq n$, with initial conditions $u(0,0)=1$ and $u(n,0)=0$ for every $n\geq 1$. Thanks to the bijection established by Callan, this recursive formula allows also to count  permutations of length $n$ avoiding $1-32-4$.  Although the recursive formula given in Equation $(\ref{(1-32-4)-recursion})$ provides an efficient way to enumerate  these permutations, we believe it provides only little insight into their structure, as it is not transparent at all from the bijection constructed by Callan how to read this recursive formula directly from a description of these permutations. Actually, it is not even clear which statistic on these permutations should correspond to the number of children of the root. Some unsuccesful attempts to read the recursive formula given in Equation $(\ref{(1-32-4)-recursion})$ directly from a description of permutations avoiding $1-32-4$ have been made. As far as we know, the best result in this direction has been achieved by Duchi, Guerrini and Rinaldi, who constructed a two labels generating tree for this class of permutations as a consequence of a certain insertion algorithm called $\mathsf{INSERTPOINT}$ (see $\cite{DGR}$). However, in the same paper, they suggest that the recursive formula given in Equation $(\ref{(1-32-4)-recursion})$ appears to be difficult to understand directly on $\mathcal{S}(1-32-4)$.  Additionally, as already pointed out by Callan,  sequence A113227 also happens to count quite a wide variety of combinatorial objects, among which we find increasing ordered rooted trees with increasing leaves, valley marked Dyck paths and inversion sequences avoiding the pattern $101$ (see $\cite{C}$ and $\cite{CMSW}$). In all these cases, it is instead relatively easy to read the recurrence relation given in Equation $(\ref{(1-32-4)-recursion})$ from the same structural description of these objects and it is actually not hard to construct quite a straightforward bijection between them. 

In this section we construct a single label generating tree for permutations avoiding the pattern $1-32-4$.  We believe that this construction finally break the annoying asimmetry between the aforementioned combinatorial objects and permutations avoiding $1-32-4$, by providing a better insight into the structure of these permutations and, compared to Callan's bijection, a clearer explanation of why the recursive formula given in Equation $(\ref{(1-32-4)-recursion})$ actually counts them. As remarkable byproducts of this construction, we also obtain an explicit algorithm to generate all permutations avoiding $1-32-4$ and we refine the enumeration of these permutations according to a simple statistic, which is the number of right-to-left maxima to the right of 1.

\begin{figure}[t]
\begin{center}\scalebox{0.5}{\begin{tikzpicture}
\draw [fill] (0,0) circle [radius=0.08];
\draw [fill] (1,1) circle [radius=0.08];
\draw [fill] (2,2) circle [radius=0.08];
\draw [fill] (3,3) circle [radius=0.08];
\draw [fill] (4,4) circle [radius=0.08];
\draw [fill] (5,3) circle [radius=0.08];
\draw [fill] (6,4) circle [radius=0.08];
\draw [fill] (7,5) circle [radius=0.08];
\draw [fill] (8,4) circle [radius=0.08];
\draw [fill] (9,3) circle [radius=0.08];
\draw [fill] (10,2) circle [radius=0.08];
\draw [fill] (11,1) circle [radius=0.08];
\draw [fill] (12,2) circle [radius=0.08];
\draw [fill] (13,1) circle [radius=0.08];
\draw [fill] (14,0) circle [radius=0.08];
\draw [fill] (15,1) circle [radius=0.08];
\draw [fill] (16,0) circle [radius=0.08];
\draw [fill,blue] (5,2) circle [radius=0.2];
\draw [fill] (5,1) circle [radius=0.08];
\draw [fill] (5,0) circle [radius=0.08];
\draw [fill,blue] (11,1) circle [radius=0.2];
\draw [fill] (11,0) circle [radius=0.08];
\draw [fill,blue] (14,0) circle [radius=0.2];
\draw[thick] (0,0)--(4,4)--(5,3)--(7,5)--(11,1)--(12,2)--(14,0)--(15,1)--(16,0);
\draw[dashed] (0,0)--(16,0);
\end{tikzpicture}}\end{center}
\caption{A valley marked Dyck path is a Dyck path in which, for each valley $DU$, one of the lattice points between the valley vertex and the $x-$axis inclusive is marked.}
\end{figure}
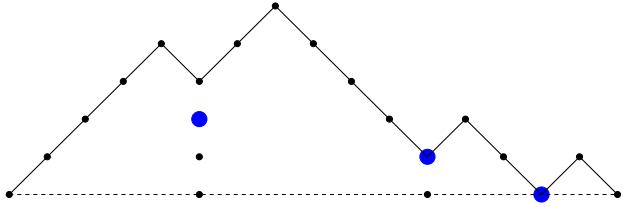

 Although not explicitly stated in $\cite{C}$, the class $\mathcal{I}$ can be described by a succession rule as follows. Consider the map $e:\mathcal{I}\longrightarrow \mathbb{N}$ attaching to every $T\in \mathcal{I}$ the number $e(T)$ of the children of its root. One can easily construct an ECO operator $\eta$ on the class $\mathcal{I}$ such that $e$ is a $\Lambda-$labelling of $\mathcal{T}_{\eta}$ where $\Lambda$ is the succession rule
$$\begin{cases} (1)\\  (k) \leadsto (1)(2)^{2}(3)^{3}...(k)^{k}(k+1)\end{cases}$$
The ECO operator $\eta$ is implicitly described by Callan in $\cite{C}$ when he proves that the class $\mathcal{I}$ is enumerated by the recursive formula in Equation $(\ref{(1-32-4)-recursion})$, but we omit further details. We will show that the class $\mathcal{S}(1-32-4)$ can be described by pretty the same rule, which  proves that the combinatorial classes $\mathcal{S}(1-32-4)$ and $\mathcal{I}$ are isomorphic. 

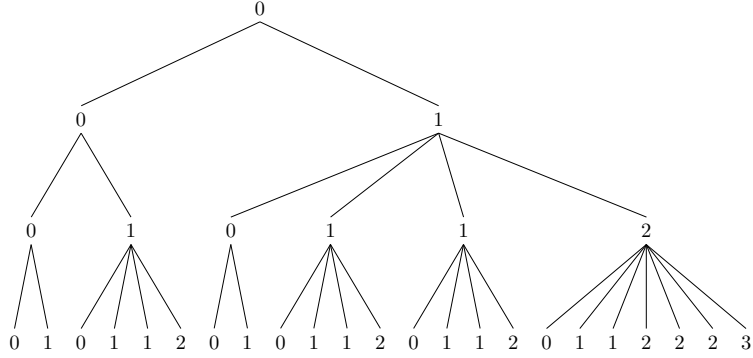
\begin{figure}\label{Gen_Tree_Rule}
\begin{center}
\scalebox{0.7}{\begin{tikzpicture}
\tikzset{level distance=60pt,sibling distance=5pt}
\Tree
[.0 [.0 [.0 [.0 ] [.1 ] ] [.1 [.0 ] [.1 ] [.1 ] [.2 ] ] ]  [.1 [.0 [.0 ] [.1 ] ] [.1 [.0 ] [.1 ] [.1 ] [.2 ] ] [.1 [.0 ] [.1 ] [.1 ] [.2 ] ] [.2 [.0 ] [.1 ] [.1 ] [.2 ] [.2 ] [.2 ] [.3 ] ] ] ]       
\end{tikzpicture}}
\end{center}
\caption{The first 3 levels of the labeled tree defined by the succession rule $\Omega$.}
\end{figure}

The outline of the proof is as follows. First, we attach to every permutation $\pi$ a label $\ell(\pi)$ defined as the number of right-to-left maxima of $\pi$ on the right of $1$ (e.g. $\ell(84617523)=3$ as the right-to-left  maxima of $84617523$ on the right of $1$ are exactly $7,5$ and $3$). Next, we  define an ECO operator $\vartheta$ on the combinatorial class $\mathcal{S}(1-32-4)$ and we show that the map $\ell:\mathcal{S}(1-32-4)\longrightarrow \mathbb{N}$ attaching to every $\pi\in \mathcal{S}(1-32-4)$ the label $\ell(\pi)$  is actually an $\Omega-$labelling of $\mathcal{T}_{\vartheta}$ where $\Omega$ is the succession rule
\begin{equation}\label{rule}
\begin{cases} (0)\\ (k) \leadsto (0)(1)^{2}(2)^{3}...(k)^{k+1}(k+1)\end{cases}
\end{equation}
This is the same as the succession rule 
$$\begin{cases} (1)\\  (h) \leadsto (1)(2)^{2}(3)^{3}...(h)^{h}(h+1)\end{cases}$$
up to the change of label $h=k+1$. The first 3 levels of the generating tree defined by  the succession rule $\Omega$ are displayed in Figure $\ref{Gen_Tree_Rule}$. Observe also that it actually takes not much effort to deduce Equation $(\ref{(1-32-4)-recursion})$ directly from  the succession rule $\Omega$, without any reference to increasing ordered rooted trees with increasing leaves. As a corollary, we find that the recursive relation in Equation $(\ref{(1-32-4)-recursion})$ for the triangle $u(n,k)$ also allows us to count permutations of length $n$ avoiding $1-32-4$ with $k-1$ right-to-left maxima to the right of 1.

Actually, it is easier to construct our ECO operator $\vartheta$ moving backwards, i.e. by first defining a reduction  operator $\rho$ on $\mathcal{S}(1-32-4)$ and then setting $\vartheta(\pi)=\rho^{\longleftarrow}(\pi)$ for every $\pi\in \mathcal{S}(1-32-4)$. 
Suppose $\pi$ is a permutation of length $n\geq 1$ avoiding $1-32-4$ such that $\ell(\pi)=k$. As already noted in $\cite{E}$, any permutation of this kind can be written in the form 
\begin{equation}\label{structure}
\pi=m_{1}\ell_{11}...\ell_{1k_{1}}m_{2}\ell_{21}...\ell_{2k_{2}} ... m_{h}\ell_{h1}...\ell_{hk_{h}}
\end{equation}
where $m_{1},...,m_{h}$ are the left to right minima of $\pi$ (and of course $m_{h}=1$) for some $h\geq 1$, while, for every $1\leq i\leq h$, the letters  $\ell_{i1},...,\ell_{ik_{i}}$ (where possibly $k_{i}=0$, with obvious meaning) denote non-empty increasing sequences such that $\max(\ell_{ij})>\max(\ell_{i(j+1)})$ for every $j\in [k_{i}-1]$ when $k_{i}\geq 2$. In particular  $k_{h}=k$ by definition of $\ell(\pi)$. Actually, a permutation that can be written in the form displayed in $(\ref{structure})$ avoids $1-32-4$ if and only if $\max(\ell_{(i+1)1})<\max(\ell_{i(k_{i}-1)})$ whenever $h\geq 2$, $i\in [h-1]$ and $k_{i}\geq 2$. For instance, the permutation $\pi=(8,9,14,12,5,2,4,10,11,1,3,13,6,7)$ avoids $1-32-4$ and decomposes as follows 
$$
\begin{matrix}
\pi=&  ( \bf{8}, & \fbox{9,14}, & \fbox{12}, & \bf{5}, & \bf{2}, & \fbox{4,10,11}, & \bf{1}, & \fbox{3,13}, & \fbox{6,7} )\\
  &  m_{1} & \ell_{11} & \ell_{12} & m_{2} & m_{3} & \ell_{31} & m_{4} & \ell_{41} & \ell_{42} 
\end{matrix}
$$
\begin{figure}\label{Pattern}
\begin{center}
\scalebox{0.3}{\begin{tikzpicture}
\draw[step=1,black,thin] (1,1) grid (14,14);
\draw [fill,red] (1,8) circle [radius=0.3];
\draw [fill] (2,9) circle [radius=0.2];
\draw [fill] (3,14) circle [radius=0.2];
\draw [fill] (4,12) circle [radius=0.2];
\draw [fill,red] (5,5) circle [radius=0.3];
\draw [fill,red] (6,2) circle [radius=0.3];
\draw [fill] (7,4) circle [radius=0.2];
\draw [fill] (8,10) circle [radius=0.2];
\draw [fill] (9,11) circle [radius=0.2];
\draw [fill,red] (10,1) circle [radius=0.3];
\draw [fill] (11,3) circle [radius=0.2];
\draw [fill] (12,13) circle [radius=0.2];
\draw [fill] (13,6) circle [radius=0.2];
\draw [fill] (14,7) circle [radius=0.2];
\end{tikzpicture}}
\end{center}
\caption{A plot of the permutation $\pi=(8,9,14,12,5,2,4,10,11,1,3,13,6,7)$ which avoids $1-32-4$. The left-to-right minima of $\pi$ are marked in bold.}
\end{figure}
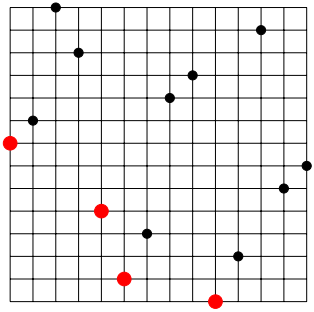
where the left-to-right minima of $\pi$ are marked in bold. Note that in this case $h=4$, while $k_{1}=2$, $k_{2}=0$, $k_{3}=1$ and $k_{4}=k=2$. 
    There is quite a natural way to reduce $\pi$ to another permutation $\rho(\pi)$ of length $n-1$ avoiding the pattern $1-32-4$, 
namely by deleting $1$, performing some further almost forced operations to restore the avoidance of the pattern $1-32-4$ and taking the standard reduction of the sequence obtained in this way (i.e. subtracting $1$ to each item of the sequence). Indeed, we can distinguish two cases:
\begin{itemize}

\item[(i)] Suppose $2$ occurs to the left of $1$ in $\pi$, i.e. $h\geq 2$ and $m_{h-1}=2$. In this case we say that $\pi$ has $\emph{type}$ $(2,1)$ and we can construct $\rho(\pi)$ as follows. We delete $1$ from $\pi$ and restore the structure displayed in $(\ref{structure})$ by sorting the list $\ell_{(h-1)1},...,\ell_{(h-1)k_{h-1}},$ $\ell_{h1},...,\ell_{hk_{h}}$ of increasing sequences to the right of $2$ in such a way that their maximum elements are in decreasing order from left to right. Finally, we define $\rho(\pi)$ as the standard reduction of the integer sequence obtained in this way. It is clear by construction that $\rho(\pi)$ avoids $1-32-4$.

\item[(ii)] Suppose $2$ occurs to the right of $1$ in $\pi$, i.e. either $h=1$ or $h\geq 2$ and $m_{h-1}\geq 3$, so that $\ell_{hi}=2\ell'_{hi}$ for some $i\in [k]$ and a possibly empty increasing sequence $\ell'_{hi}$. In this case we say that $\pi$ has $\emph{type}$ $(1,2)$ and  we can construct $\rho(\pi)$ as follows. We delete $1$ and restore the structure displayed in $(\ref{structure})$ by moving $2$ to the position previously occupied by $1$ (in other words we swap $1$ and $2$ and then delete $1$), thus obtaining the integer sequence $m_{1}\ell_{11}...\ell_{1k_{1}}m_{2}\ell_{21}...\ell_{2k_{2}} ... 2\ell_{h1}...\ell_{hi}'...\ell_{hk_{h}}$. Finally, we define $\rho(\pi)$ as the standard reduction of this sequence. Again, it is clear by construction that $\rho(\pi)$ avoids $1-32-4$. 

\end{itemize}

\begin{example} Let us illustrate the previous construction with two examples. 
\begin{itemize}
\item[(i)] Take the permutation $\pi=(8,9,14,12,5,2,4,10,11,1,3,13,6,7)$ and let us compute its reduction $\rho(\pi)$. Note that $\pi$ has type $(2,1)$, therefore we first delete $1$ to obtain the sequence $(8,9,14,12,5,2,4,10,11,3,13,6,7)$, then we sort the increasing sequences $(4,10,11)$, $(3,13)$ and $(6,7)$ on the right of $2$ in such a way that their maximum elements $11,13$ and $7$ are in decreasing order. Therefore the correct order is given by $(3,13)(4,10,11)(6,7)$, which yields the sequence $(8,9,14,12,5,2,3,13,4,10,11,6,7)$. Taking the standard reduction of this sequence returns the permutation $\rho(\pi)=(7,8,13,11,4,1,2,12,3,9,10,5,6)$. 
\item[(ii)] Take the permutation $\pi=(8,9,14,12,5,3,4,10,11,1,6,13,2,7)$ and let us compute its reduction $\rho(\pi)$. Note that $\pi$ has type $(1,2)$, hence we first delete $1$ to obtain the sequence $(8,9,14,12,5,3,4,10,11,6,13,2,7)$, then we move $2$ to the position previously occupied by $1$, so to obtain the sequence $(8,9,14,12,5,3,4,10,11,2,6,13,7)$. Taking the standard reduction of this sequence returns the permutation $\rho(\pi)=(7,8,13,11,4,2,3,9,10,1,5,12,6)$. 
\end{itemize}
\end{example}

This construction induces a reduction operator $\rho$ on $\mathcal{S}(1-32-4)$ and thus, as mentioned before, an ECO operator $\vartheta$ on $\mathcal{S}(1-32-4)$. Now we want to show that $\vartheta$ can be described by the succession rule $\Omega$ given by ($\ref{rule}$). For this purpose, we will explicitly describe all the elements of $\vartheta(\pi)$ and compute their labels by reversing the previous construction. More specifically, we will expand $\pi$ by appending $0$ at the end of $\pi$, then moving some of the increasing sequences $\ell_{h1},...,\ell_{hk_{h}}$ to the right of $0$ and finally normalizing the sequence thus obtained (i.e. adding $1$ to each item of the sequence). In fact, the range of possibilities to perform this operation is quite constrained because we have to preserve the avoidance of the pattern $1-32-4$. First append a $0$ at the end of $\pi$.
\begin{itemize}
\item[(i)] Of course, no occurrence of $1-32-4$ will appear if we move the whole sequence $\ell_{h1}...\ell_{hk_{h}}$ to the right of $0$. In this way we get the sequence $m_{1}\ell_{11}...\ell_{1k_{1}}...m_{h}0\ell_{h1}...\ell_{hk_{h}}$, whose normalization is a permutation which we denote by $\pi^{(k)}$. Note that $\ell(\pi^{(k)})=k$.
\item[(ii)] Suppose instead that $k\geq 1$ and we want to move only $i\in \{0,...,k-1\}$ increasing sequences among $\ell_{h1},...,\ell_{hk_{h}}$ to the right of $0$.  Then it is easy to see that there is a unique way to perform this operation so to prevent an occurrence of $1-32-4$ to appear in the resulting  expansion of $\pi$, which is the following way.  Choose some $j\in [i+1]$ and move the $(i+1)^{th}$ suffix $\ell_{h(k_{h}-i)},...,\ell_{hk_{h}}$ of the list $\ell_{h1},...,\ell_{hk_{h}}$, except for its $j^{th}$ increasing sequence $\ell_{h(k_{h}-i+j-1)}$, to the right of $0$.  In other words, move the sequence $\ell_{h(k_{h}-i)}...\hat{\ell}_{h(k_{h}-i+j-1)}...\ell_{hk_{h}}$ (where the hat over an item means that it must be omitted) to the right of $0$, thus obtaining the sequence $$m_{1}\ell_{11}...\ell_{1k_{1}}...m_{h}\ell_{h1}...\ell_{h(k_{h}-i-1)}\ell_{h(k_{h}-i+j-1)}0\ell_{h(k_{h}-i)}...\hat{\ell}_{h(k_{h}-i+j-1)}...\ell_{hk_{h}}.$$ Finally normalize this sequence to a permutation, which we denote by $\pi^{(i,j)}$. Note that $\ell(\pi^{(i,j)})=i$.  Hence, this operation produces $i+1$ children with label $i$, for every $0\leq i\leq k-1$, from a node with label $k$. 
\end{itemize}

Note that all permutations defined in $(i)$ and $(ii)$ will have type $(2,1)$, therefore these permutations cannot exhaust the whole class $\mathcal{S}(1-32-4)$ and we need to construct other expansions of $\pi$ to include also permutations of type $(1,2)$. To this purpose we also move $m_{h}=1$ to the right of $0$ and perform some further transformations.  

\begin{itemize}
\item[(iii)] Of course, no occurrence of $1-32-4$ will appear if we move the whole sequence $1\ell_{h1}...\ell_{hk_{h}}$ to the right of $0$. In this way, we get the sequence $m_{1}\ell_{11}...\ell_{1k_{1}}...01\ell_{h1}...\ell_{hk_{h}}$, whose normalization is a permutation which we denote by $\pi^{[1]}$. More generally, it is clear that no occurrence of $1-32-4$ will appear if we perform the following operation. Choose some $i\in [k]$ and move $1$ immediately to the left of $\ell_{hi}$, then move the sequence $\ell_{h1}...1\ell_{hi}...\ell_{hk_{h}}$ to the right of $0$. In this way we obtain the sequence $m_{1}\ell_{11}...\ell_{1k_{1}}...0\ell_{h1}...1\ell_{hi}...\ell_{hk_{h}}$, whose normalization is a permutation, which we denote by $\pi^{[i]}$. Note that $\ell(\pi')=k$. Hence, operation $(i)$ and $(iii)$ produce $k+1$ children with label $k$ from a node with label $k$.
\item[(iv)] Finally, we have a last possibility to transform $\pi$ and prevent an occurrence of $1-32-4$ to appear. Move $1$ back to the right of $\ell_{hk_{h}}$, then move the sequence $\ell_{h1}...\ell_{hk_{h}}1$ to the right of $0$. In this way we obtain the sequence $m_{1}\ell_{11}...\ell_{1k_{1}}...0\ell_{h1}...\ell_{hk_{h}}1$ and normalize it to a permutation, which we denote by $\pi^{[k+1]}$.  Note that in this case, unlike  in the previous one, we have $\ell(\pi^{[k+1]})=k+1$. Hence, this operation will produces a unique child with label $k+1$ from a node with label $k$. Note that this last possibility  could actually be regarded as a special case of $(iii)$ if we let $\pi$ terminate with an additional empty increasing sequence $\ell_{(h+1) k_{h+1}}$.
\end{itemize}

Note that all permutations defined in $(iii)$ and $(iv)$ have type $(1,2)$.

\begin{example} Let us take $\pi=({\bf{5}},9,14,10,12,{\bf{1}},2,7,13,6,11,3,8,4)$ as a working example to illustrate some of the previous constructions. Note that in this case $\pi$ has the form $m_{1}\ell_{11}\ell_{12}m_{2}\ell_{21}\ell_{22}\ell_{23}\ell_{24}$ where $m_{1}=5$, $\ell_{12}=(9,14)$ and $\ell_{12}=(10,12)$, while $m_{2}=1$, $\ell_{21}=(2,7,13)$, $\ell_{22}=(6,11)$, $\ell_{23}=(3,8)$ and $\ell_{24}=(4)$, in particular $\ell(\pi)=4$. First we insert a $0$ at the end of $\pi$ to obtain the sequence $(5,9,14,10,12,1,2,7,13,6,11,3,8,4,0)$.
\begin{itemize}
\item[(i)] We start by constructing $\pi^{(4)}$. We move the sequence $(2,7,13,6,11,3,8,4)$ to the right of $0$, thus obtaining the sequence $(5,9,14,10,12,1,0,2,7,13,6,11,3,8,4)$, whose normalization is given by the permutation $\pi^{(4)}=(6,10,15,11,13,2,1,3,8,14,7,12,4,9,5)$.
\item[(ii)] Now let us construct the permutation $\pi^{(2,2)}$. We move the suffix 
$(6,11)(3,8)(4)$ of the list 
$(2,7,13)(6,11)(3,8)(4)$, except for its $2^{nd}$ element 
$(3,8)$, to the right of $0$, thus obtaining the sequence $(5,9,14,10,12,1,2,7,13,3,8,0,6,11,4)$, whose normalization is given by the permutation $\pi^{(2,2)}=(6,10,15,11,13,2,3,8,14,4,9,1,7,12,5)$.   
\item[(iii)] Let us now construct the permutation $\pi^{[3]}$. We move the sequence $(1,2,7,13,6,11,3,8,4)$  to the right of $0$ and move $1$ immediately to the left of $(3,8)$, thus obtaining the sequence $(5,9,14,10,12,0,2,7,13,6,11,1,3,8,4)$, whose normalization is given by the permutation  $\pi^{[3]}=(6,10,15,11,13,1,3,8,14,7,12,2,4,9,5)$. 
\item[(iv)] Finally we construct the permutation $\pi^{[4]}$. We move the sequence $(1,2,7,13,6,11,3,8,4)$  to the right of $0$ and move $1$ immediately to the right of $(4)$, thus obtaining the sequence $(5,9,14,10,12,0,2,7,13,6,11,3,8,4,1)$, whose normalization is given by the permutation $\pi^{[3]}=(6,10,15,11,13,,1,3,8,14,7,12,4,9,5,2)$. 
\end{itemize}  
\end{example}

\begin{figure}\label{Gen_Tree}
\begin{center}
\scalebox{0.5}{\begin{tikzpicture}[grow=right]
\tikzset{level distance=150pt,sibling distance=0.1pt}
\Tree
[.1(0) [.21(0) [.321(0) [.4321(0) ] [.4312(1) ] ] [.312(1) [.4231(0) ] [.4213(1) ] [.4123(1) ] [.4132(2) ]] ]  [.12(0) [.231(0) [.3421(0) ] [.3412(1) ] ] [.213(1) [.3241(0) ] [.3214(1) ] [.3124(1) ] [.3142(2) ] ] [.123(1) [.2341(0) ] [.2134(1) ] [.1234(1) ] [.1342(2) ] ] [.132(2) [.2431(0) ] [.2413(1) ] [.2314(1) ] [.2143(2) ] [.1243(2) ] [.1423(2) ] [.1432(3) ] ] ] ]       
\end{tikzpicture}}
\end{center}
\caption{The first 3 levels of the generating tree for permutations avoiding the pattern $1-32-4$, where each node $\pi\in \mathcal{S}(1-32-4)$ is labeled by $(\ell(\pi))$.}
\end{figure}
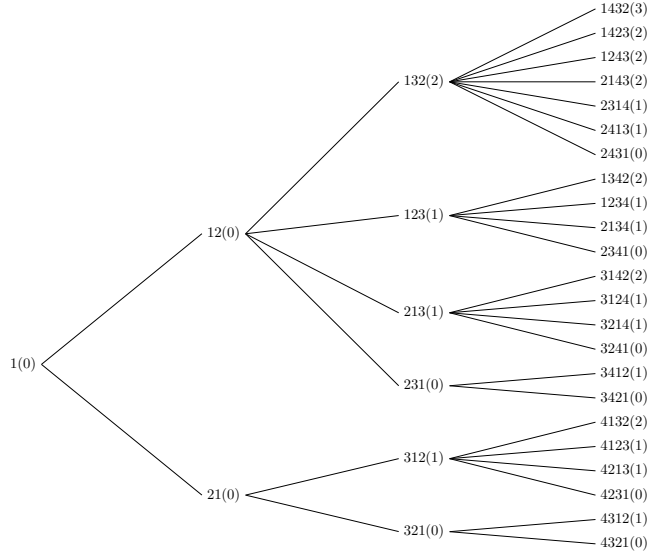

 Now we are in a position to state and prove the main result. 
 
\begin{theorem}
Suppose $\pi\in \mathcal{S}(1-32-4)$ and $k=\ell(\pi)$. 
\begin{itemize}
\item[(i)] If $k=0$, then $\vartheta(\pi)=\{\pi^{(0)},\pi^{[1]}\}$. 
\item[(ii)] If $k\geq 1$, then $\vartheta(\pi)=\{\pi^{(k)},\pi^{(i,j)},\pi^{[p]}:\ 0\leq i\leq k-1, 1\leq j\leq i+1,1\leq p\leq k+1\}$.
\item[(iii)] The map $\ell$ is an $\Omega-$labelling of $\mathcal{T}_{\vartheta}$.
\end{itemize}
 
\end{theorem} 

\proof  
It is mere routine to check that $\rho(\pi^{(0)})=\rho(\pi^{[1]})=\pi$ when $k=0$ and that $\rho(\pi^{(k)})=\rho(\pi^{(i,j)})=\rho(\pi^{[p]})=\pi$ when $k\geq 1$,  $0\leq i\leq k-1, 1\leq j\leq i+1 $ and $1\leq p\leq k+1$. Conversely, assume $\sigma\in \vartheta(\pi)$ and write $\sigma$ in the form $m_{1}\ell_{11}...\ell_{1k_{1}}...m_{h}\ell_{h1}...\ell_{hk_{h}}$ as in Equation $(\ref{structure})$. Suppose first that $\sigma$ has type $(2,1)$. If $2$ and $1$ occur consecutively in $\pi$, then it is immediate to see that $\sigma=\pi^{(k)}$. Otherwise, it is also fairly easy to check that $\sigma=\pi^{(i,j)}$ where $i=\ell(\sigma)$, while $j=1$ in case $\max(\ell_{h1})<\max(\ell_{(h-1)k_{h-1}})$ and $j=\max\{j\in \{2,3,...,i+1\}:\ \max(\ell_{h(j-1)})>\max(\ell_{(h-1)k_{h-1}})\}$ otherwise. Suppose now $\sigma$ has type $(1,2)$. Then it is also easy to check that $\sigma=\pi^{[p]}$ where $p$ is the unique $p\in [k_{h}]$ such that $2$ occurs in the increasing sequence $\ell_{hp}$. This proves $(i)$ and $(ii)$. Finally, $(iii)$ holds because we know that $\ell(\pi^{(0)})=0$ and $\ell(\pi^{[1]})=1$ when $k=0$, while $\ell(\pi^{(i,j)})=i$, $\ell(\pi^{(k)})=\ell(\pi^{[p]})=k$ and $\ell(\pi^{[k+1]})=k+1$ when $k\geq 1$,  $0\leq i\leq k-1, 1\leq j\leq i+1 $ and $1\leq p\leq k$.
\endproof
 
We observe once again that the construction above actually provides an algorithm to generate all permutations avoiding $1-32-4$. The permutations generated this way up to length $4$ are displayed in Figure $\ref{Gen_Tree}$. We end this section with the following corollary summarizing the further byproduct of the previous construction, which is the refined enumeration of permutations avoiding $1-32-4$ according to the number of right-to-left maxima to the right of 1.

\begin{corollary} For every $0\leq k\leq n-1$ denote by $v(n,k)$ the number of permutations avoiding $1-32-4$ with length $n$ and $k$ right-to-left maxima to the right of $1$. Then the triangle $v(n,k)$ satisfies the recurrence relation  
$$v(n,k)=v(n-1,k-1)+(k+1)\sum_{j=k}^{n-2}v(n-1,j)$$
where we agree that $v(0,-1)=1$ and $v(n,-1)=0$ for every $n\geq 1$. 
\end{corollary}

\begin{table}
\begin{center}
\begin{tabular}{|c|cccccccc|}
  \hline
  $n\backslash k$ & 0 & 1 & 2 & 3 & 4 & 5 & 6 & 7  \\
  \hline
  1 & 1 & & & & & & & \\
  2 & 1 & 1 & & & & & & \\ 
  3 & 2 & 3 & 1 & & & & & \\
  4 & 6 & 10 & 6 & 1 & & & & \\ 
  5 & 23 & 40 & 31 & 10 & 1 & & & \\ 
  6 & 105 & 187 & 166 & 75 & 15 & 1 & & \\
  7 & 549 & 993 & 958 & 530 & 155 & 21 & 1 & \\
  8 & 3207 & 5865 & 5988 & 3786 & 1415 & 287 & 28 & 1\\ 
  \hline
\end{tabular}
\end{center}
\caption{The triangle $v(n,k)$ counting permutations avoiding $1-32-4$ with length $n$ and  $k$ right-to-left maxima to the right of $1$, for $1\leq n\leq 8$ and $0\leq k\leq n-1$.}
\label{m1}
\end{table}


\section{Conclusion and further work}

In this paper we finally succeded in describing a generating tree with single label for permutations avoiding $1-32-4$. However, as pointed out in Section $\ref{1-32-4avoiders}$, the sequence counting these permutations also happens to count other combinatorial structures, such as increasing ordered rooted trees with increasinig leaves or valley marked Dyck paths, and  describing a clear bijection between these structures would have a great combinatorial significance. Although a bijection was already established by Callan in $\cite{C}$ through a sequence of non-trivial steps and despite the remarkable advances in $\cite{DGR}$ to make more explicit the connection between these objects, we believe it is worth looking for a bijection admitting a reasonably simpler description. Using the generating tree described in this paper, it is likely that one can recursively construct some other bijection, preserving the respective labels, and we hope this could eventually lead to a significant advance in this direction. 

As for a more generic issue, we observe that permutations avoiding the classical pattern $1324$ form a subclass of permutations avoiding the vincular pattern $1-32-4$.  Hence, it might be interesting to investigate in which cases the ECO operator $\vartheta$ fails to expand a permutation avoiding the classical pattern $1324$ to another permutation avoiding the same pattern, causing an occurrence of $1324$ to appear. 

Furthermore, we believe it is worth investigating whether the construction provided in this paper to generate all permutations avoiding $1-32-4$ can be adapted to generate all permutations avoiding the vincular pattern $1-(\pi+1)-(|\pi|+2)$, for some particular consecutive patterns $\pi$ other than $12$ or $21$ (where $\pi+1$ denotes the permutation obtained from $\pi$ by adding $1$ to each item of $\pi$). In this case there would be hope to find a recurrence relation to count these permutations. 

 In general,  the fast recurrence for permutations avoiding $1-32-4$ provided in $\cite{C}$ suggests to look for nonobvious recurrences counting other similar patterns. For instance, note that a kind of structural description for permutations avoiding $12-34$ is given in $\cite{E}$, hence there may be hope to find a nice generating tree and deduce some reasonable recurrence relation, just like we did in the case of permutations avoiding $1-23-4$.  

Finally, still concerning enumeration of permutations avoiding $1-32-4$, although the recursive formula found by Callan provides quite an efficient and elegant way to count them, a closed formula would clearly be a more satisfactory answer. In this regard, the generating function $u(z)$ of permutations avoiding $1-32-4$ can be recursively described as a continued fraction $u(z)=1-z(U(0)-z)$ where $U(n)=1-z^{n}-z/U(n+1)$ for every $n\in\mathbb{N}$ (see $\cite{S}$). However, this description can hardly be considered a closed formula. The succession rule describing permutations avoiding $1-32-4$ can also be translated into a functional equation for the generating function $u(z,t)$  of these permutations (where $z$ keeps track of the length and $t$ keeps track of the label), as explained at the very end of Section $\ref{ECO}$. This equation is actually a linear PDE of the form
\begin{equation}\label{PDE}
(1-t)zt^{2}\frac{\partial u}{\partial t}(z,t)+((1-t)^{2}(1-zt)+zt)u(z,t)=zt(1-t)^{2}+ztu(z,1).
\end{equation}
However, we do not know whether Equation $(\ref{PDE})$ provides enough information to find some closed form expression for $u(z,t)$. Further research in this direction could improve our understanding of sequence A113227.

\paragraph*{Acknowledgement} The author whishes to thank Luca Ferrari for kindly reviewing some technical details related to the construction illustrated in this paper and for his valuable suggestions to improve the organization of its content.

\end {document}